\documentclass[12pt]{article}
\usepackage{amsxtra}
\usepackage{amssymb, amsmath}
\usepackage{amsthm}
\usepackage {hyperref}

\usepackage{amssymb}
\usepackage{amsmath}
\usepackage[dvips]{graphicx}

\begin{document}
\begin{center}
\textbf{\LARGE{Upside Down Numerical Equation, Bimagic Squares, and the  day  September 11}}
\end{center}

\bigskip
\begin{center}
\textbf{\large{Inder Jeet Taneja}}\\
Departamento de Matem\'{a}tica\\
Universidade Federal de Santa Catarina\\
88.040-900 Florian\'{o}polis, SC, Brazil.\\
\textit{e-mail: ijtaneja@gmail.com\\
http://www.mtm.ufsc.br/$\sim$taneja}
\end{center}

\bigskip
\begin{abstract}
\textit{In this short note we have given an equation based on the date 11.09.2001 and presented some magic squares. The magic squares presented are of order 3  $ \times $ 3, 4  $\times $ 4, 5  $\times $ 5, 9  $ \times $ 9 , 16   $\times $ 16 and 25  $\times $ 25. While the magic square of higher order 9  $ \times $ 9 , 16   $ \times $ 16 and 25  $\times $ 25 are bimagic. These magic squares are based on the digits, 1-6-9, 0-1-2, 0-1-2-9 and 0-1-2-6-9. The interesting fact in all these magic squares is that if we gave them a rotation of 180$^{o}$, they remain again the magic squares. In order to have this rotation, we have used the numbers in digital forms. Moreover, the day of submission of this work (20.10.2010) has only the digits 0, 1 and 2. Using only these three digits, we have presented magic squares of order 3 $\times$ 3 and 9 $\times$ 9.}
\end{abstract}

\section{Numerical Equation}

Here below we shall present a numerical equation based on the dates 9-11-2001 to 9-11-2010. We generally write, either 9/11 or 11/9 . Total period from 2001 to 2010 is of 9 years.

{Let us write it as:  \hspace*{25pt} \includegraphics[bb=0mm 0mm 208mm 296mm, width=6.6mm, height=3.9mm, viewport=3mm 4mm 205mm 292mm]{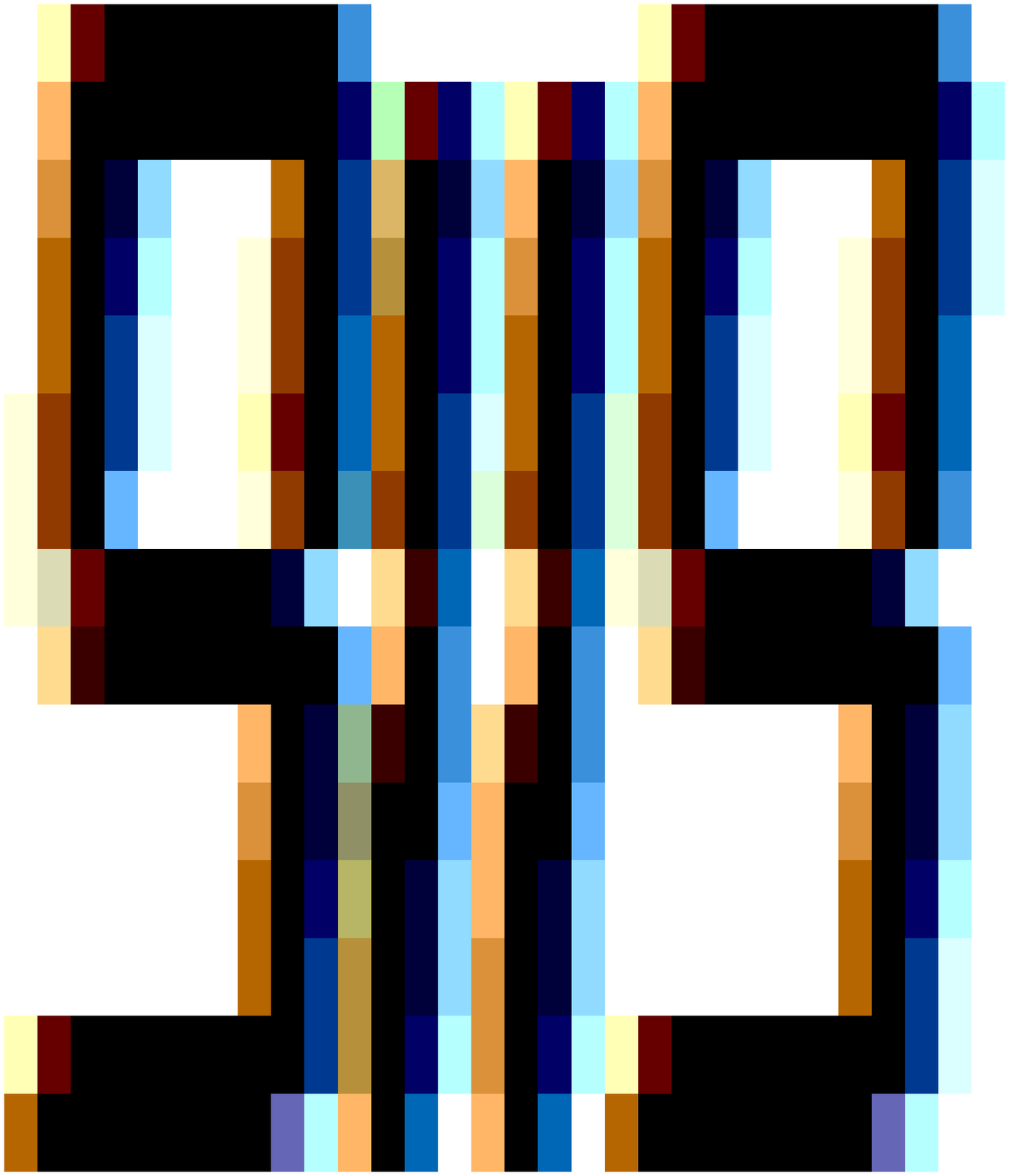}

$180^{o}$ (degrees) rotation:  \includegraphics[bb=0mm 0mm 208mm 296mm, width=6.4mm, height=4.3mm, viewport=3mm 4mm 205mm 292mm]{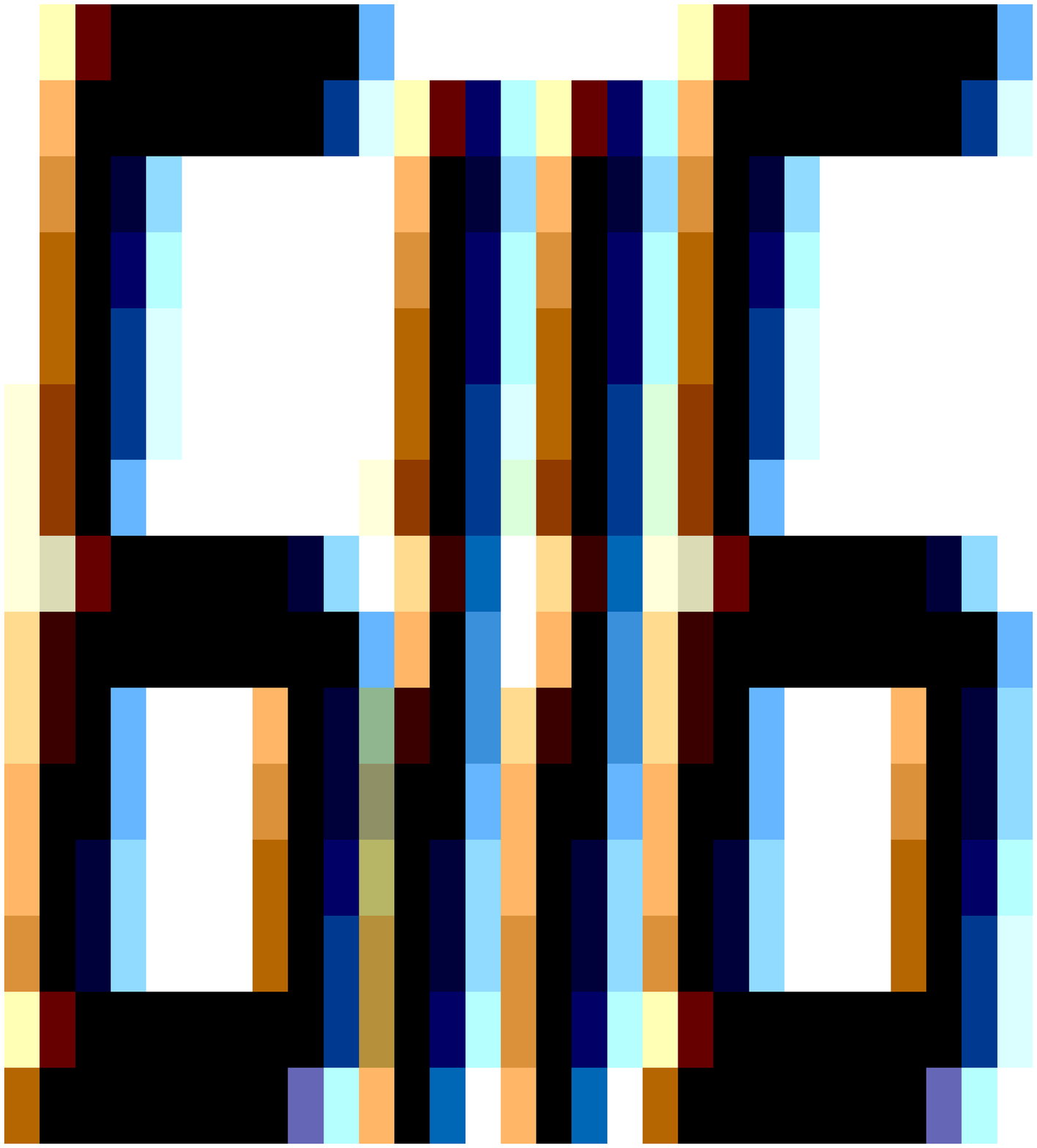}

Year:     \hspace*{90pt}  \includegraphics[bb=0mm 0mm 208mm 296mm, width=7.8mm, height=4.3mm, viewport=3mm 4mm 205mm 292mm]{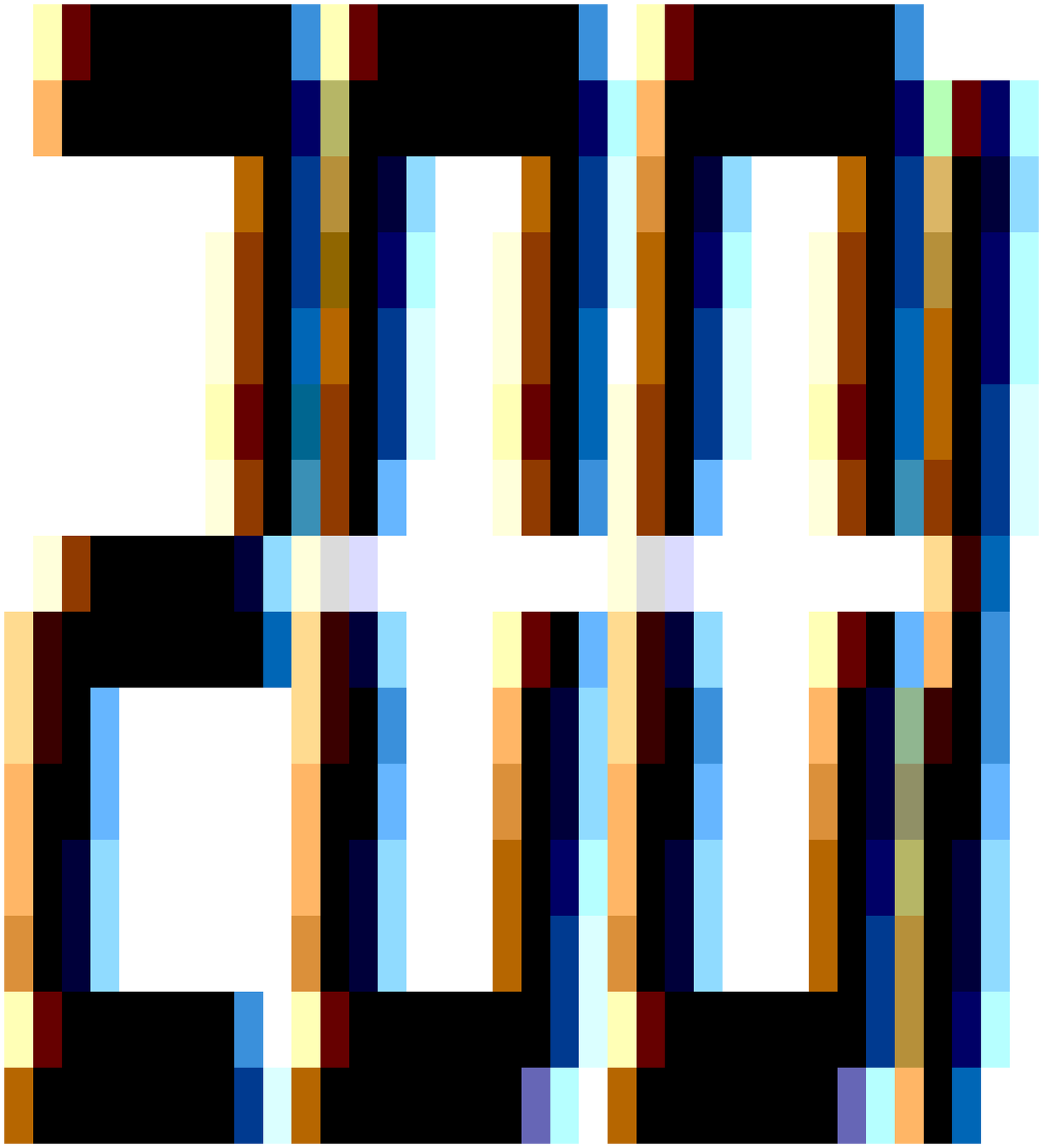}

$180^{o}$ (degrees) rotation:  \includegraphics[bb=0mm 0mm 208mm 296mm, width=7.8mm, height=3.9mm, viewport=3mm 4mm 205mm 292mm]{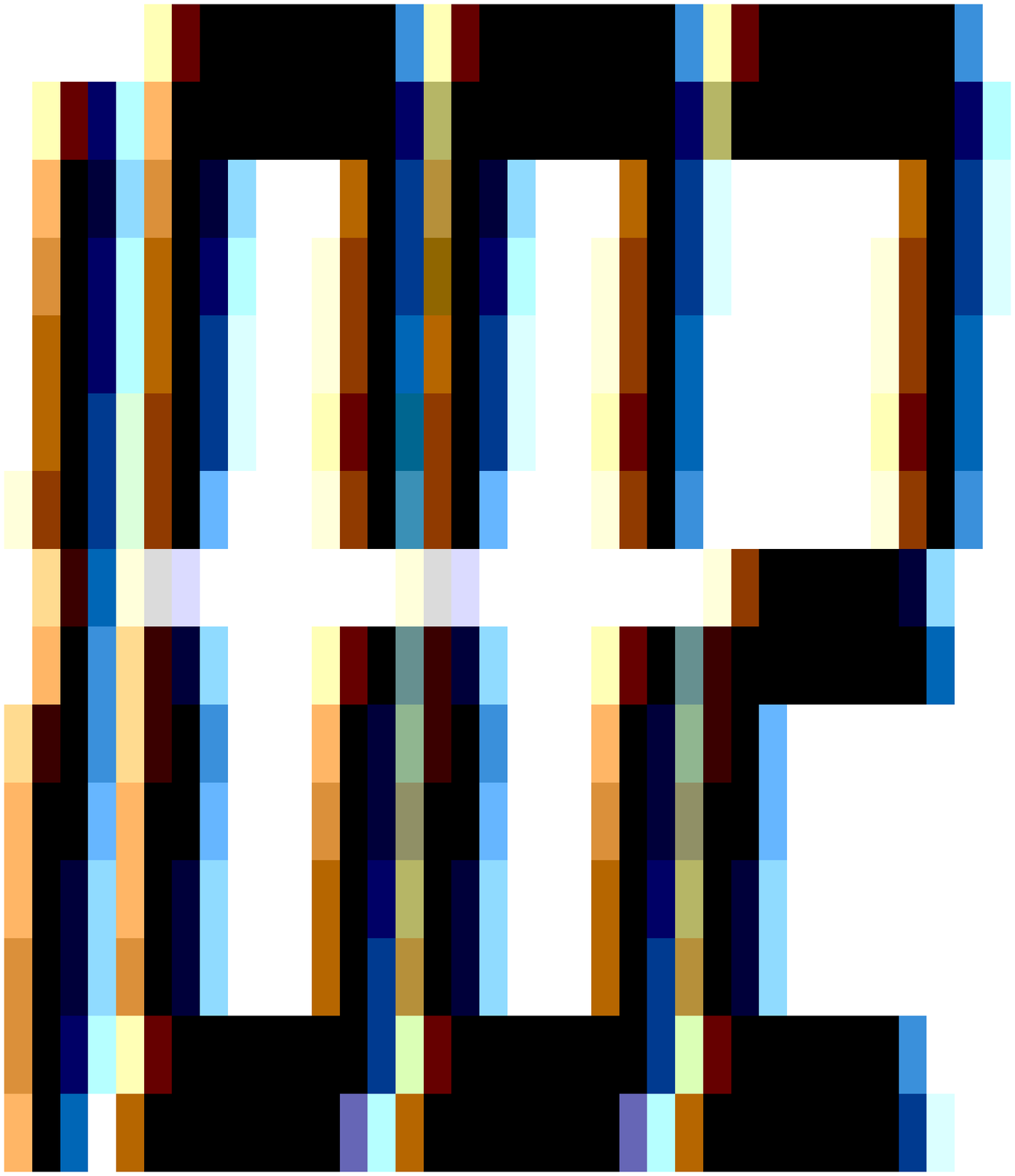}.}

\bigskip
Based on above four numbers, we have the following upside down numerical equation:
\begin{center}
\textbf{\textit{\includegraphics[bb=0mm 0mm 208mm 296mm, width=43.8mm, height=4.5mm, viewport=3mm 4mm 205mm 292mm]{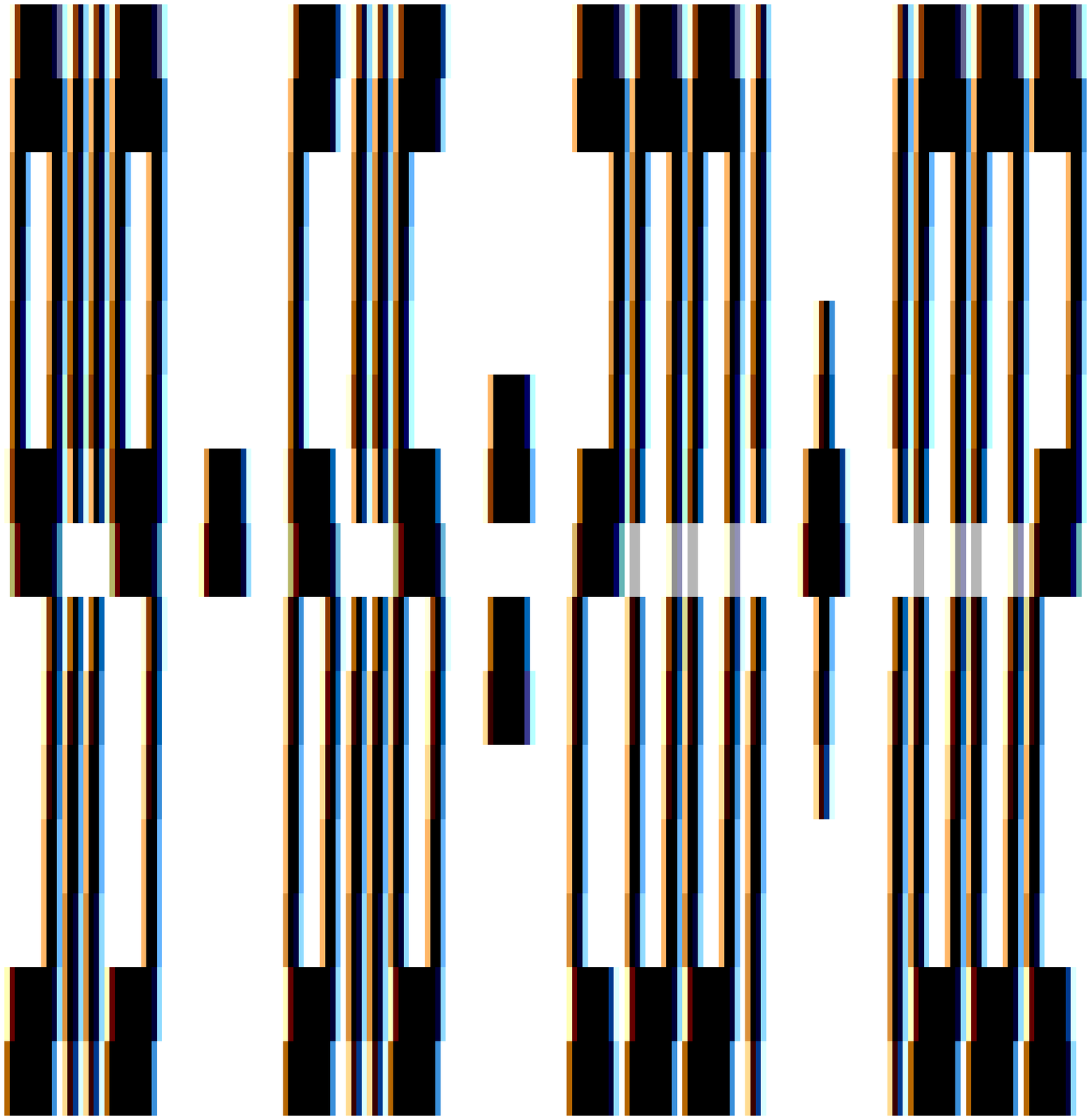}}}\textit{.}\textbf{\textit{}}
\end{center}

Here we have used the numbers in the digital forms. This allows us to make rotation of $180^{o}$ and interestingly the equation remains the same.

The digits appearing on the right side of above equation 0, 1 and 2 are the same on day 20.10.2010, i.e., the day of submission of  this work.

Now based on the digits appearing in the above equation we shall give below magic squares, while some of them are bimagic.

\subsection{Universal, Upside Down and Bimagic Squares}

\bigskip
Here below are definitions of Universal and Bimagic squares.
\bigskip

\noindent \textbf{\large{$\bullet$ Bimagic square}}

\begin{enumerate}
\item [(i)] A \textbf{magic square} is a collection of numbers put as a square matrix, where the sum of elements of each row, sum of elements of each column and sum of elements of each of two principal diagonals is always same.  For simplicity, let us write it as $\bf{S1}$.

\item  [(ii)] \textbf{bimagic square} is a magic square where the sum of square of each element of rows, columns and two principal diagonals are the same. For simplicity, let us write it as $\bf{S2}$.
\end{enumerate}

\bigskip
\noindent \textbf{\large{$\bullet$ Universal magic square}}

\bigskip
Universal magic squares are the magic squares having following properties:

\begin{enumerate}
\item [(i)] \textbf{Upside Down}, i.e. if we rotate it $180^{o}$ (degrees) it remains again the magic square;

\item [(ii)] \textbf{Mirror looking}, i.e., if we put it in front of mirror or see from the other side of the glass, or see on the other side of the paper, it always remains the magic square.
\end{enumerate}

\section{Upside Down Magic Squares of Order $\bf{3 \times 3}$ and $\bf{9 \times 9}$}

In the above equation we have the digits \textbf{1-6-9} on the left side and  the digits \textbf{0-1-2} on the right side. In both cases we have below \textit{upside down magic squares}  formed by two algorism combinations made from three digits in each case. See below:

\bigskip
\noindent \textbf{$\bullet$ Semi-magic square of order $ \bf{3 \times 3} $ having only the digits 1,6 and 9}

\begin{center}
\includegraphics[bb=0mm 0mm 208mm 296mm, width=45.8mm, height=35.6mm, viewport=3mm 4mm 205mm 292mm]{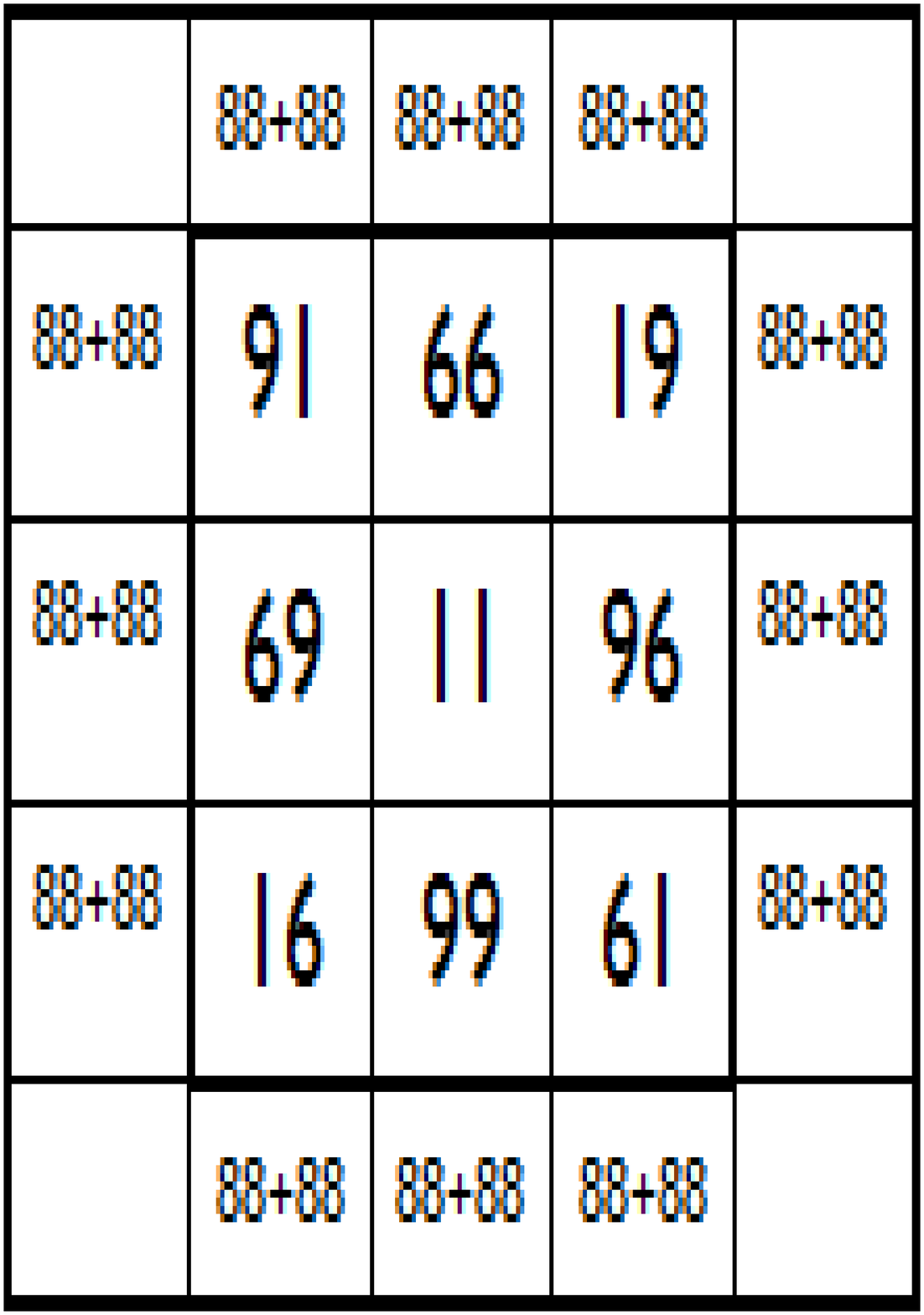}
\end{center}

\newpage
\noindent \textbf{$\bullet$ Magic square of order $ \bf{3 \times 3} $ having only the digits 0,1 and 2}

\begin{center}
\includegraphics[bb=0mm 0mm 208mm 296mm, width=46.6mm, height=32.4mm, viewport=3mm 4mm 205mm 292mm]{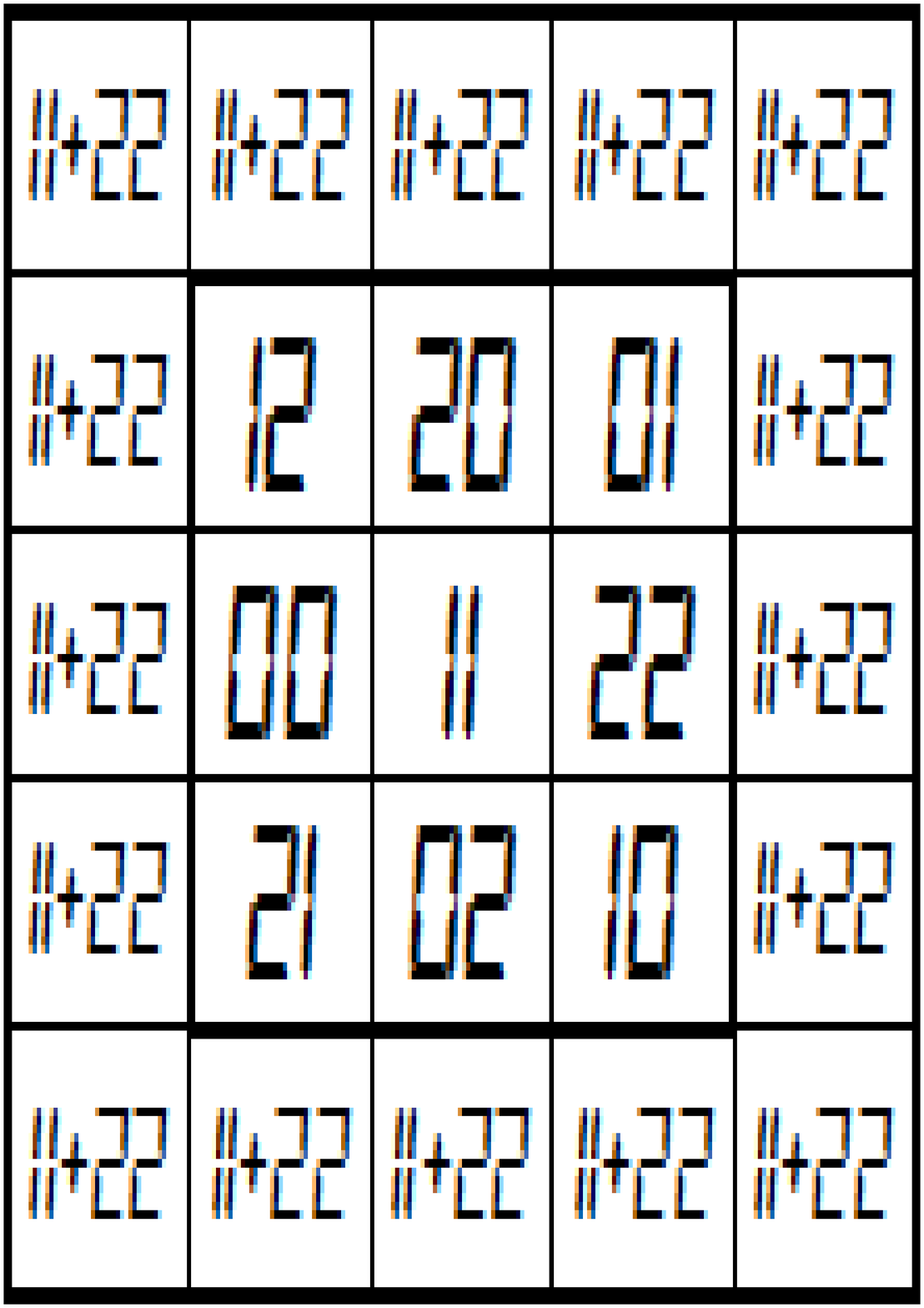}
\end{center}

\bigskip
We note that the magic square made from the digits 1, 6 and 9 is not complete, i.e., it is \textbf{\textit{semi-magic square}} since the sum of the principal diagonals is not the same. While, the second one is a \textbf{\textit{complete magic square}}. Here the sums in both the cases are also \textit{upside down}, i.e., they remains the same when we give a rotation of $180^{o}$. Moreover, the second magic square can be considered as universal, i.e., when we look it from the mirror it is again a magic square, obviously, in this case the sum is not the same. More studied on \textit{universal magic squares} can be seen in Taneja \cite{tan3}.

\subsection{Upside Down Bimagic Squares of Order $ \bf{9 \times 9}$}

The above two magic squares of order $ 3\times 3$ are made only from two algorism combinations of 1-6-9 and 0-1-2 respectively. Instead, if we consider four algorism combinations of these digits, we can make $3 \times 3 \times 3 \times 3= 9 \times 9= 81$ different numbers in each case. Interesting these numbers lead us to the following two \textit{upside down bimagic squares}.

\bigskip
\noindent \textbf{$\bullet$ Upside down bimagic square of order $ \bf{9 \times 9} $ using only the digits 1, 6 and 9}
\bigskip

\begin{center}
\textbf{\includegraphics[bb=0mm 0mm 208mm 296mm, width=122.2mm, height=48.7mm, viewport=3mm 4mm 205mm 292mm]{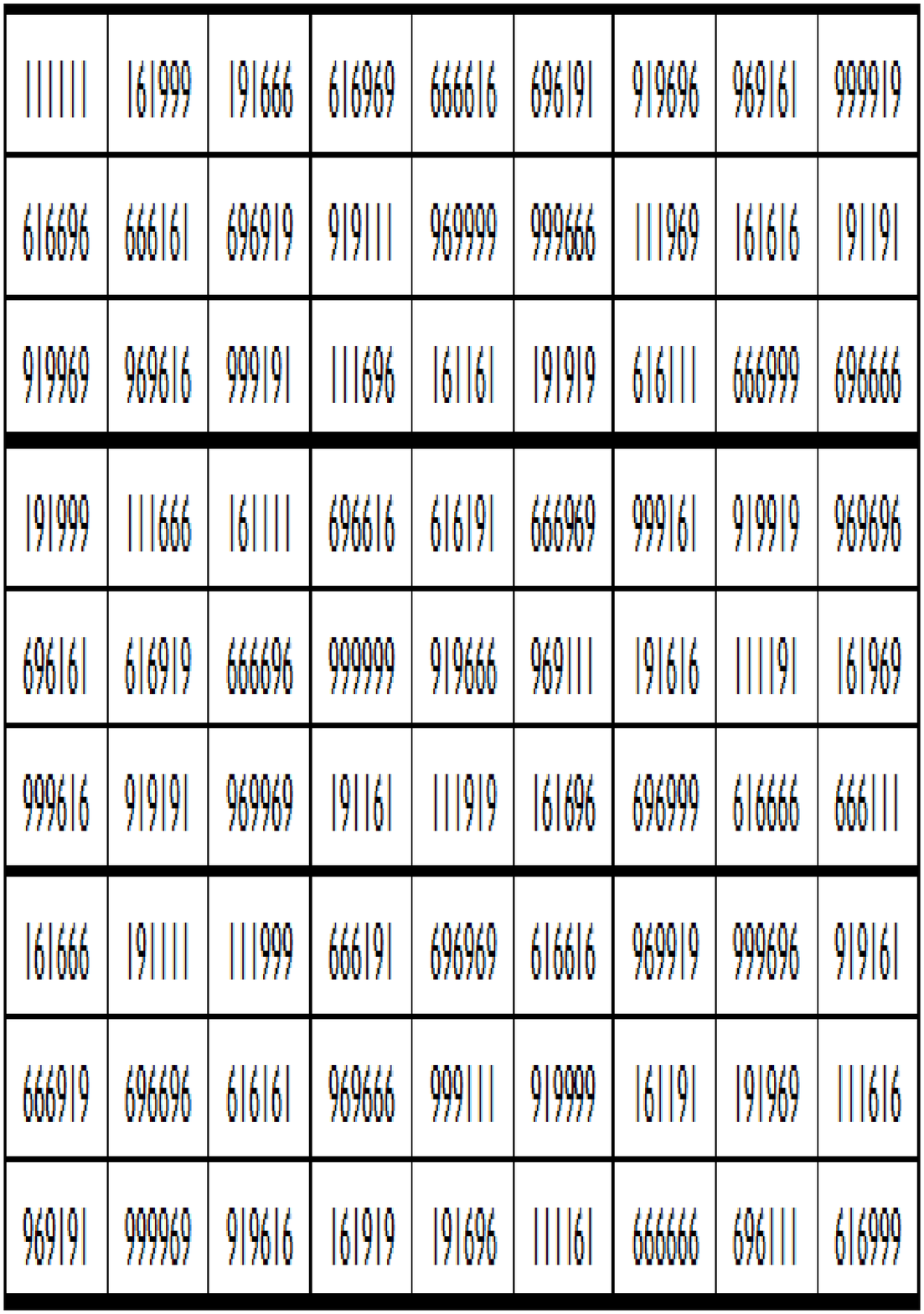}}
\end{center}

Here we have  $S1 := 53328$ and $S2 := 414976074$. Also we have sum of each block of $3 \times 3=53328$ and the sum of square of each term in each block of  $3 \times 3=414976074$.

\newpage
\noindent \textbf{$\bullet$ Upside down bimagic square of order $ \bf{9 \times 9} $ using only the digits 0, 1 and 2}
\bigskip

\begin{center}
\includegraphics[bb=0mm 0mm 208mm 296mm, width=127.4mm, height=51.0mm, viewport=3mm 4mm 205mm 292mm]{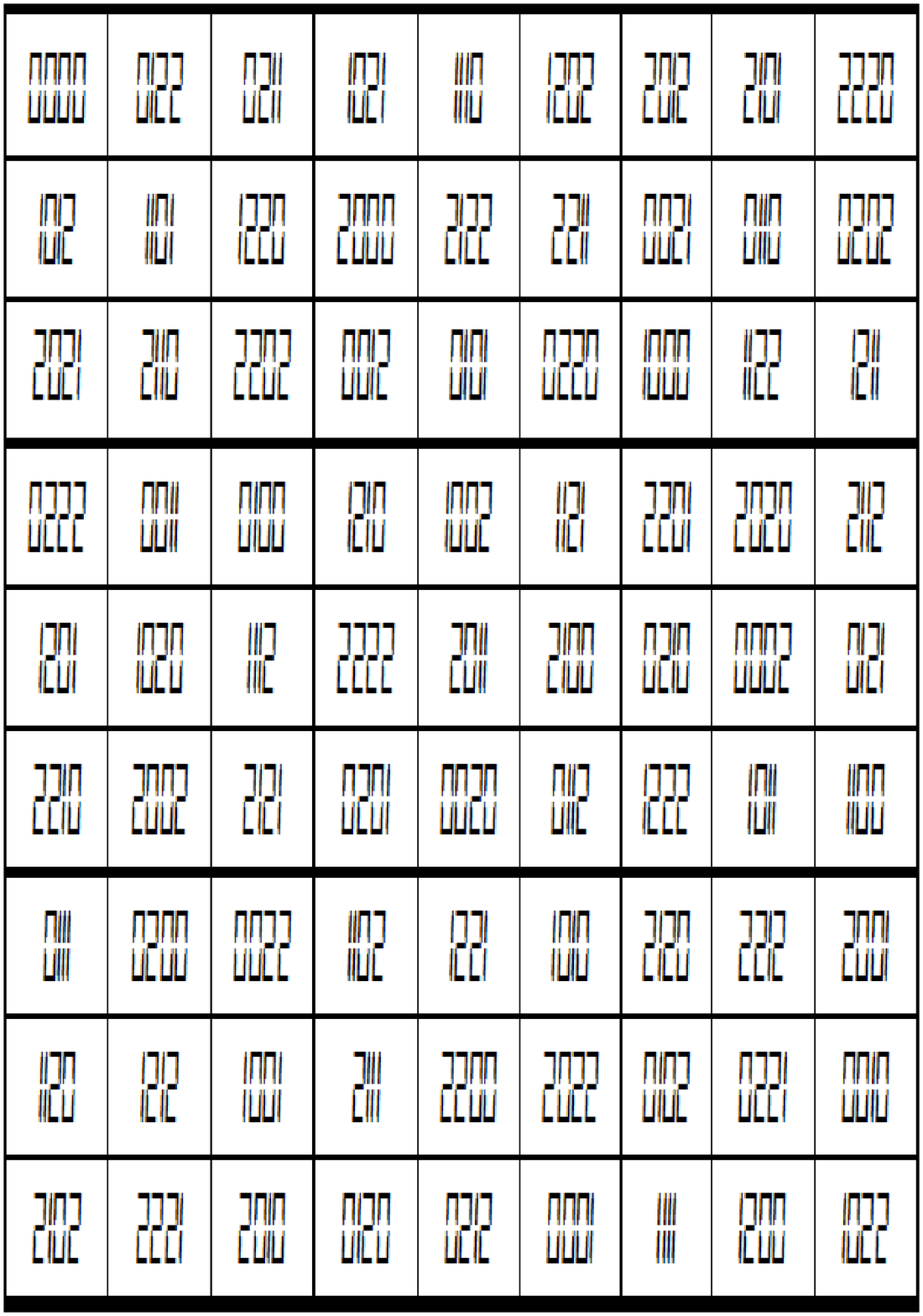}
\end{center}

Here we have $S1 := 9999$ and $S2 := 17169395$. Also we have sum of each block of  $3 \times 3=9999$ and the sum of square of each term in each block of  $3 \times 3=17169395$.

\section{Upside Down Bimagic Square of Order $\bf{16 \times 16}$}

This section deals with the bimagic square of order $16 \times 16 $ based on the digits appearing in the day 11.01.2001. A magic square of order $4 \times 4$ is also presented.

\bigskip
\noindent \textbf{\large{$\bullet$ The dates 11.09.2001 or 11.09.2010}}
\bigskip

We observe that the above dates are formed by four digit, i.e.,
\begin{center}
\textbf{0, 1, 2 and 9}
\end{center}

Making two algorism combination of above four digits, we have 16 different numbers. These 16 numbers give us the following upside down magic square of order $ 4 \times 4$ of sum $S1:=132$. See below:

\begin{center}
\includegraphics[bb=0mm 0mm 208mm 296mm, width=32.8mm, height=27.1mm, viewport=3mm 4mm 205mm 292mm]{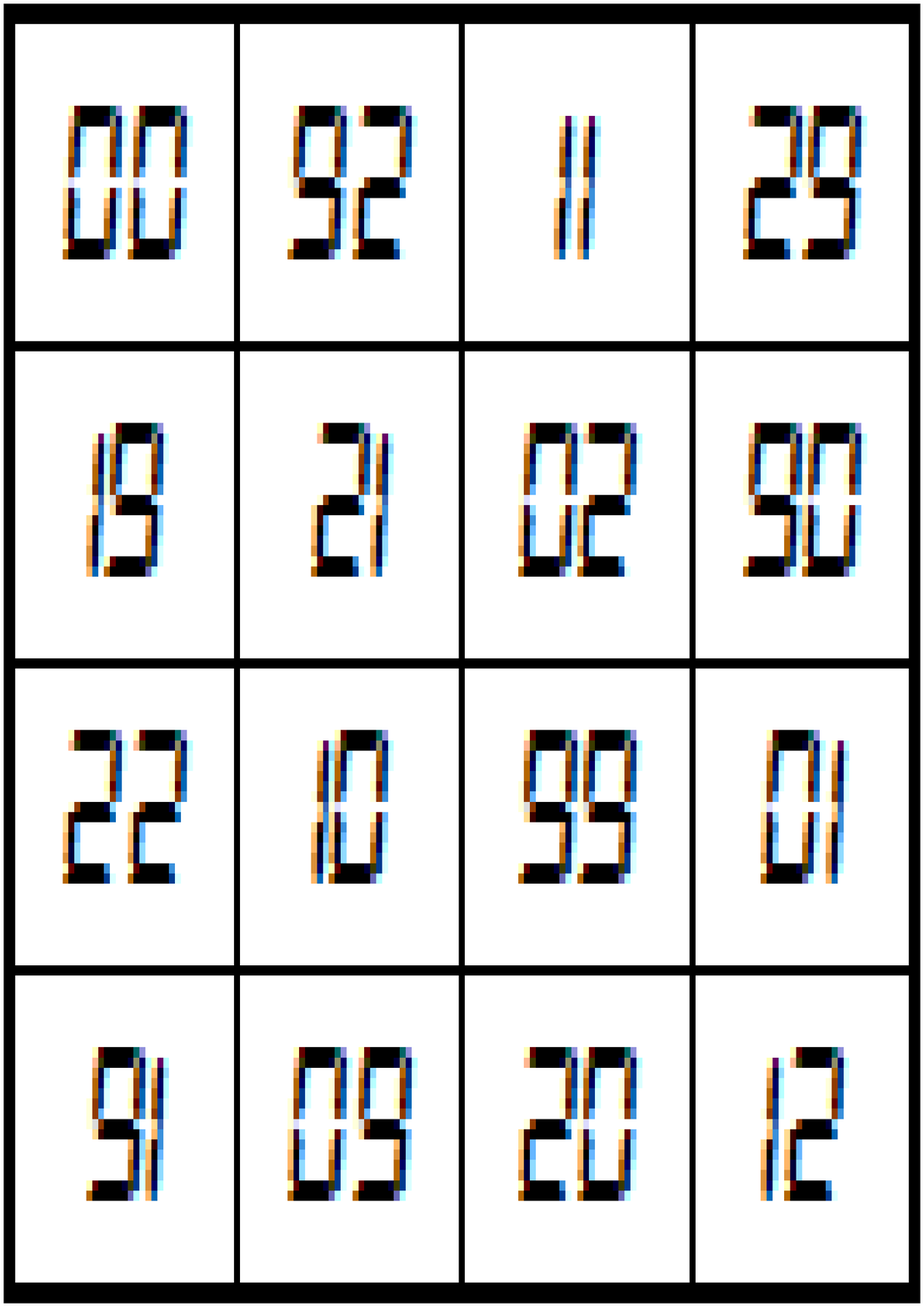}
\end{center}

In we make four algorism combination from the four digits, 0, 1, 2 and 9, we can make $4 \times 4 \times 4 \times 4 = 256$ different numbers and these 256 different numbers. Also we can write $ 256 = 16 \times 16$. Based on this we can construct an \textbf{\textit{ upside down bimagic square}} of order $ 16 \times 16$. See below:

\begin{center}
\includegraphics[bb=0mm 0mm 208mm 296mm, width=148.8mm, height=76.6mm, viewport=3mm 4mm 205mm 292mm]{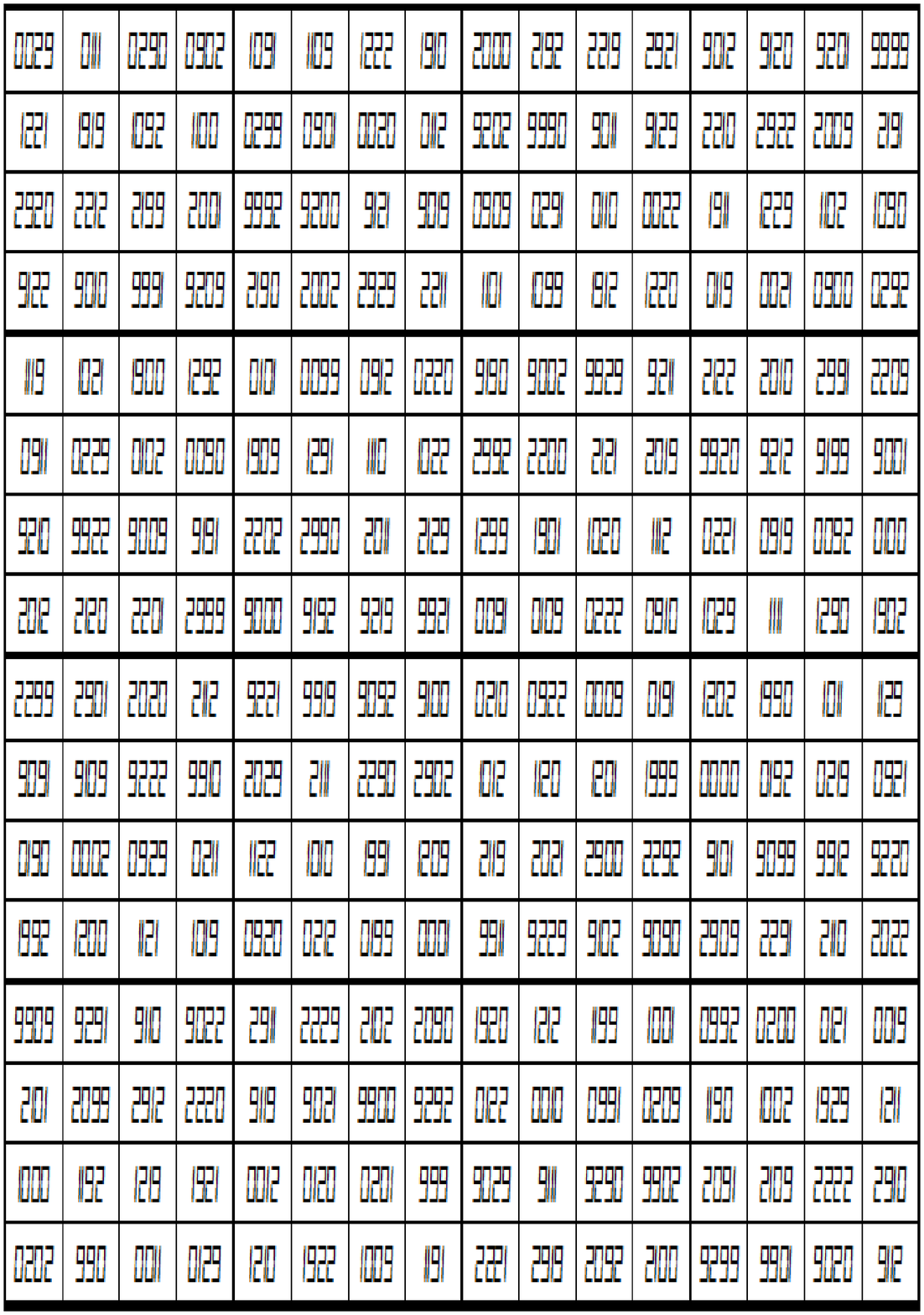}
\end{center}

Here we have $ S1 := 53328 $ and  $ S2 := 379762424 $. Also we have sum of each block of order $4 \times 4 = 53328 $ and the sum of square of each term in each block of $ 4 \times 4 = 379762424 $.

\bigskip
\textbf{Remark 1.} It is interesting to note that the sum $S1_{9\times 9} $ made from the digits 1, 6 and 9 and the sum $S1_{16\times 16} $ made from the digits 0, 1, 2 and 9 is the same, i.e.,

\[S1_{16\times 16} =S1_{9\times 9} = 53328 .\]

\section{Upside Down Bimagic Square of Order $ \bf{25 \times 25}$}

We observe that the numerical equation given in section 1 has only the following five digits:

\begin{center}
\textbf{0, 1, 2, 6 and 9.}
\end{center}

Here below is an upside down magic square of order $ 5 \times 5$ of sum $S1 := 198$ made from two algorism combination of above five digits.

\begin{center}
\includegraphics[bb=0mm 0mm 208mm 296mm, width=37.9mm, height=38.1mm, viewport=3mm 4mm 205mm 292mm]{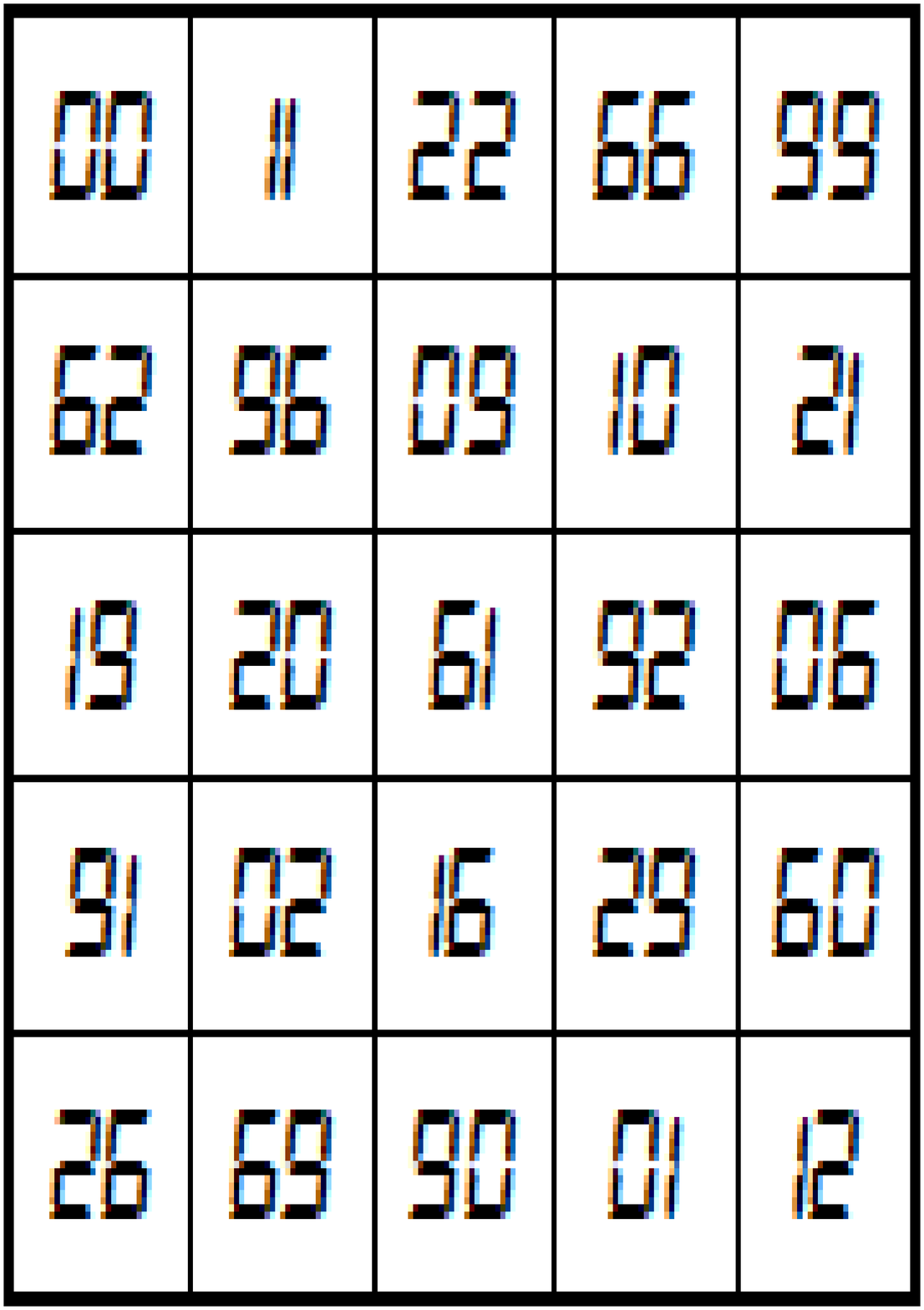}
\end{center}

In we make four algorism combination of above five digits, 0, 1, 2, 6 and 9 we can have $ 5 \times 5 \times 5 \times 5 = 625$ different numbers and these 625 different numbers can be written as $625 = 25 \times 25$. Based on this we can construct an \textbf{\textit{ upside down bimagic square}} of order $ 25 \times 25$. See below:
\begin{center}
\includegraphics[bb=0mm 0mm 208mm 296mm, width=160.6mm, height=90.5mm, viewport=3mm 4mm 205mm 292mm]{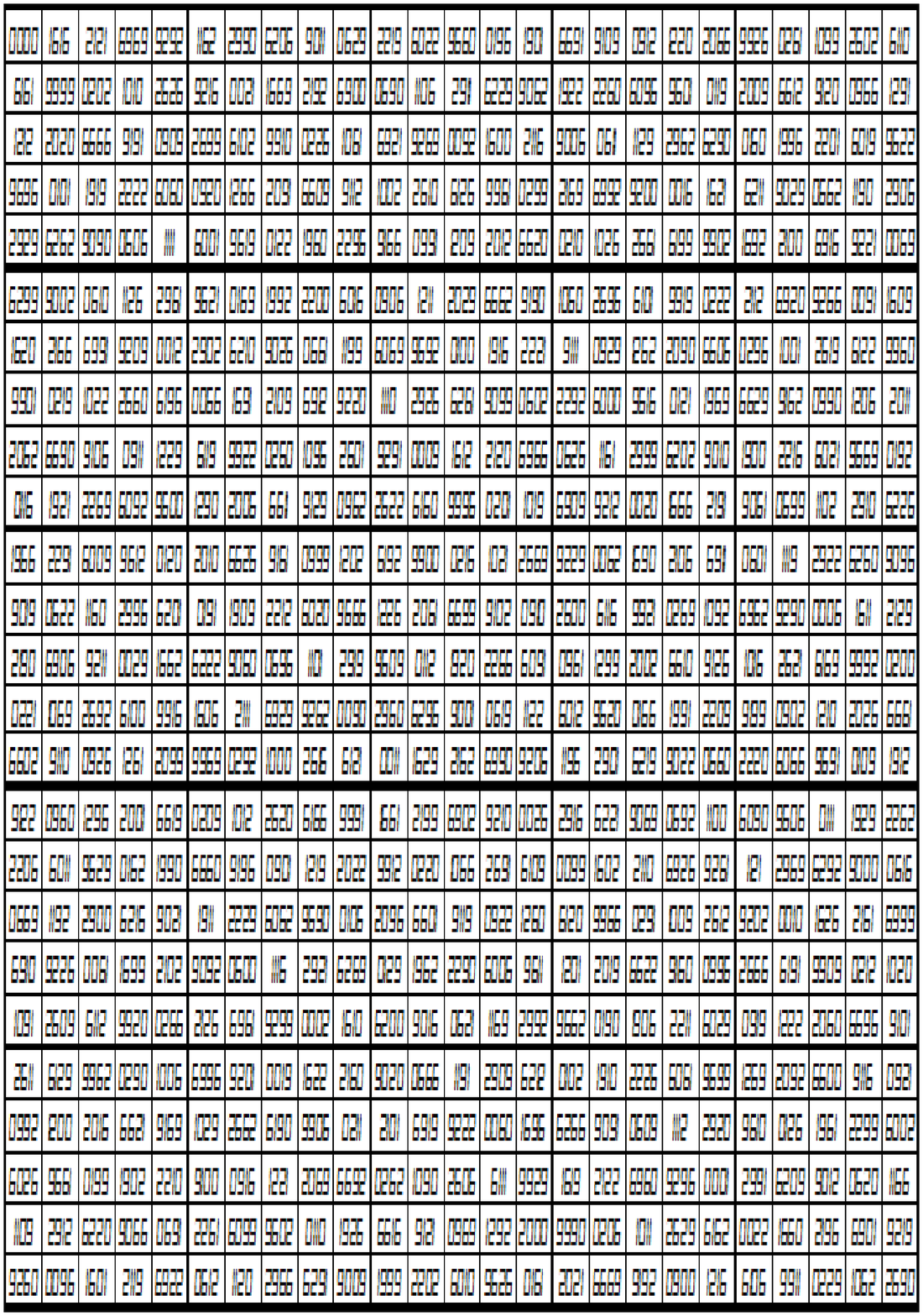}
\end{center}

Here we have $S1 := 99990$ and $S2 := 688808890$. Also we have sum of each block of order $5 \times 5 = 99990$ and the sum of square of each term in each block of order $5 \times 5 = 688808890$.

\bigskip
We observe that the above \textit{upside down magic square} of order $25 \times 25$ is also is \textit{pandiagonal}.

\bigskip
\textbf{Remark 2. } It is interesting to note that the sum $S1_{9\times 9} :=9999$ of the magic square of order $9 \times 9$ made from the digits 0, 1 and 2 and $S1_{25\times 25} :=99990$ of the magic square of order $25 \times 25$ made from five digits 0,1,2, 6 and 9 has the following interesting relation:

\[S1_{25\times 25} =10\times S1_{9\times 9} =99990.\]

\section{Final Comments}

We have considered here the dates 11.09.2001 or 11.09.2010 just as natural curiosity. There are many other dates that also have the same digits, such as, 09.11.2010, 09.01.2010, 19.02.2010, etc. There is a big collection of work on \textit{bimagic squares}  in the site maintained by C. Boyer \cite{boyer}. The another site on magic squares maintained by H. Heinz \cite{hei}  also brings a good collection of curiosities, books, sites, papers etc. We have brought here for the first in the literature an idea of \textit{upside down bimagic square of higher order}. Some recent work by author on \textit{universal bimagic squares} are given in Taneja \cite{tan4}. If we consider seven digits 0, 1, 2, 5, 6, 8 and 9 in the digital form we can bring an \textit{upside down bimagic square of order $ 49 \times 49$}. This study is given in a separate work by author \cite{tan5}. Some interesting relations, curiosities and magic squares with the digits 1, 6 and 9 are under preparation Taneja \cite{tan6}.

\begin{center}
---------------------------
\end{center}

\end{document}